\documentclass[10pt,twocolumn,superscriptaddress,nofootinbib]{revtex4}
\usepackage[T1]{fontenc}
\usepackage[utf8]{inputenc}
\usepackage{color}
\usepackage{amssymb}
\usepackage{booktabs}
\usepackage{mathrsfs}
\usepackage{url}
\usepackage{hyphenat}
\usepackage{mathtools}
\usepackage{cancel} 
\usepackage{ifpdf}
\usepackage{mdwtab}
\usepackage{mymathabx}
\usepackage{enumitem}
\ifpdf
  \usepackage[pdftex]{graphicx}
  \usepackage[pdftex,margin=0.6in,letterpaper]{geometry}
  \usepackage[bookmarks=true, bookmarksopen=true,%
    bookmarksdepth=3,bookmarksopenlevel=2,%
    colorlinks=true,%
    linkcolor=blue,%
    citecolor=blue,%
    filecolor=blue,%
    menucolor=blue,%
    urlcolor=blue]{hyperref}
  \hypersetup{pdftitle={A tour of bordered Floer theory}}
  \hypersetup{pdfauthor={Robert Lipshitz, Peter S. Ozsváth, and Dylan
      P. Thurston}}
\else
  \usepackage[dvips]{graphicx}
  \usepackage[draft]{hyperref}
\fi
\usepackage{tikz}

\usetikzlibrary{matrix,arrows}
\tikzstyle{every picture}=[> = latex']
\tikzset{cdlabel/.style={above,sloped,
    execute at begin node=$\scriptstyle,execute at end node=$}}
\tikzset{algarrow/.style={->, thick}}
\tikzset{blgarrow/.style={->, thick}}
\tikzset{clgarrow/.style={->, thick}}
\tikzset{tensoralgarrow/.style={->, thin, double}}
\tikzset{tensorblgarrow/.style={->, thin, double}}
\tikzset{tensorclgarrow/.style={->, thin, double}}
\tikzset{modarrow/.style={->, dashed}}
\tikzset{othmodarrow/.style={->, thick}}
\tikzset{Amodar/.style={->, dashed}}
\tikzset{Dmodar/.style={->, dashed}}

\newcommand{\RR}{\mathbb R}

\newcommand{\ZZ}{\mathbb Z}

\newcommand{\FF}{\mathbb F}

\newcommand{\bD}{\mathbb{D}}

\newcommand{\co}{:}

\newcommand{\OneHalf}{{\textstyle\frac{1}{2}}}

\newcommand{\bdy}{\partial}

\newcommand{\lbracket}{[}
\newcommand{\rbracket}{]}

\newcommand{\spinc}{\mathfrak s}

\DeclareMathOperator{\Sym}{Sym}

\DeclareMathOperator{\Hom}{Hom}

\DeclareMathOperator{\Tor}{Tor}

\DeclareMathOperator{\spin}{spin}
\newcommand{\SpinC}{\spin^c}

\DeclareMathOperator{\inv}{inv}
\DeclareMathOperator{\Inv}{Inv}
\DeclareMathOperator{\ev}{ev}
\DeclareMathOperator{\gr}{gr}

\DeclareMathOperator{\Cone}{Cone}

\newtheorem{theorem}{Theorem}

\newtheorem{example}[theorem]{Example}
\newtheorem{lemma}[theorem]{Lemma}
\newtheorem{remark}[theorem]{Remark}
\newtheorem{definition}[theorem]{Definition}

\newcommand{\HF}{\mathit{HF}}
\newcommand{\HFa}{\widehat {\HF}}

\newcommand{\CF}{{\mathit{CF}}}
\newcommand{\CFa}{\widehat {\mathit{CF}}}
\newcommand{\HFKa}{\widehat{\mathit{HFK}}}
\newcommand{\x}{\mathbf x}
\newcommand{\y}{\mathbf y}

\newcommand{\w}{\mathbf w}

\newcommand\HH{\mathit{HH}}

\newcommand\Hochschild\HH

\newcommand{\Ainf}{\mathcal A_\infty}
\newcommand{\Alg}{\mathcal{A}}

\newcommand{\Idem}{\mathcal{I}}
\newcommand{\alphas}{{\boldsymbol{\alpha}}}
\newcommand{\betas}{{\boldsymbol{\beta}}}

\newcommand{\rhos}{{\boldsymbol{\rho}}}

\newcommand{\bSigma}{\widebar{\Sigma}}

\newcommand{\cM}{\mathcal{M}}

\newcommand{\tcM}{\widetilde{\mathcal{M}}}
\newcommand{\DD}{\textit{DD}}
\newcommand{\DA}{\textit{DA}}

\newcommand{\CFD}{\mathit{CFD}}
\newcommand{\CFDD}{\mathit{CFDD}}
\newcommand{\CFA}{\mathit{CFA}}

\newcommand{\CFDA}{\mathit{CFDA}}
\newcommand{\CFDAa}{\widehat{\CFDA}}
\newcommand{\CFAA}{\mathit{CFAA}}
\newcommand{\CFAAa}{\widehat{\CFAA}}
\newcommand{\CFDa}{\widehat{\CFD}}

\newcommand{\CFDDa}{\widehat{\CFDD}}

\newcommand{\CFAa}{\widehat{\CFA}}

\newcommand{\cZ}{\mathcal{Z}}
\newcommand{\PtdMatchCirc}{\cZ}
\newcommand{\PMC}{\PtdMatchCirc}

\newcommand{\CircPts}{{\mathbf{a}}}

\newcommand{\dg}{\textit{dg} }

\newcommand\Id{\mathbb{I}}
\newcommand\Ground{\mathbf k}

\newcommand\DTP{\mathop{\widetilde\otimes}\nolimits}
\newcommand\DT{\boxtimes}
\newcommand\Gen{\mathfrak{S}}

\newcommand{\Field}{{\FF_2}}

\newcommand{\Heegaard}{\mathcal{H}}
\newcommand{\HD}{\Heegaard}

\newcommand{\HB}{\mathsf{H}}

\newcommand{\smallGroup}{G}

\DeclareMathOperator{\Mor}{Mor}

\newcommand{\op}{\mathrm{op}}

\newcommand\rKh{\widetilde{\mathit{Kh}}}
\newcommand\Kh{\mathit{Kh}}

\makeatletter
\newcommand\honestalg[3]{\bigl\lbracket
\begin{smallmatrix} #1\@ifempty{#3}{}{&#3} \\ #2 \end{smallmatrix}
\bigr\rbracket}

\makeatother

\newcommand{\sos}[3]{\mathbin{{}_{#1}\mathord#2_{#3}}}
\newcommand{\lsub}[2]{{}_{#1}#2}

\newcommand{\smargin}[1]{%
  \marginpar[\raggedleft\footnotesize{#1}]{\raggedright\footnotesize{#1}}%
}

\newread\testin

\graphicspath{{draws/}{mpdraws/}{}}
\makeatletter
\def\input@path{{}{draws/}}
\makeatother

\def\mathcenter#1{%
  \vcenter{\hbox{$#1$}}%
}

\DeclareRobustCommand{\widebar}[1]{\overline{#1}{}}

\makeatletter
\newcommand\mi@kern[1]{%
  \settowidth\@tempdima{$\mi@obj^{#1}$}
  \kern-\@tempdima
  #1
  \settowidth\@tempdima{$\mi@obj$}
  \kern\@tempdima
}

\newtoks\mi@toksp
\newtoks\mi@toksb
\DeclareRobustCommand{\manyindices}[5]{
  \def\mi@obj{#5}
  \mi@toksp\expandafter{\mi@kern{#2}}
  \mi@toksb\expandafter{\mi@kern{#1}}
  \@mathmeasure4\textstyle{#5_{#1}^{#2}}
  \@mathmeasure6\textstyle{#5_{#3}^{#4}}
  \dimen0-\wd6 \advance\dimen0\wd4
  \@mathmeasure8\textstyle{\hphantom{{}_{#1}^{#2}}#5^{\the\mi@toksp#4}_{\the\mi@toksb#3}}
  \hbox to \dimen0{}{\kern-\dimen0\box8}
}
\makeatother 

\usepackage{ifpdf}
\ifpdf
  
\else
  
\fi

\newcommand{\HMto}{\widecheck{\mathit{HM}}}
\newcommand{\HMfrom}{\widehat{\mathit{HM}}}
\newcommand{\HMbar}{\widebar{\mathit{HM}}}
\renewcommand{\smargin}[1]{} 
\newenvironment{proofsketch}{\textit{Proof sketch.}}{\hfill$\Box$}

\begin{document}
\title{A tour of bordered Floer theory}

\author{Robert Lipshitz}
\thanks{RL was supported by NSF grant DMS-0905796, the Mathematical
  Sciences Research Institute, and a Sloan Research
  Fellowship}
\email{lipshitz@math.columbia.edu}
\affiliation{Department of Mathematics, Columbia University\\
  New York, NY 10027}

\author{Peter~S.~Ozsv\'ath}
\thanks{PSO was supported by NSF grant DMS-0505811, the Mathematical
  Sciences Research Institute, and a Clay
  Senior Scholar Fellowship}
\email {petero@math.mit.edu}
\affiliation {Department of Mathematics, MIT\\ Cambridge, MA 02139}

\author{Dylan~P.~Thurston}
\thanks {DPT was supported by NSF
  grant DMS-1008049 and the Mathematical
  Sciences Research Institute}
\email{dthurston@barnard.edu}
\affiliation{Department of Mathematics,
         Barnard College,
         Columbia University\\
         New York, NY 10027}

\begin{abstract} Heegaard Floer theory is a kind of topological
    quantum field theory, assigning graded groups to closed,
    connected, oriented 3-manifolds and group homomorphisms to
    smooth, oriented 4-dimensional cobordisms. Bordered Heegaard
    Floer homology is an extension of Heegaard Floer homology to
    3-manifolds with boundary, with extended-TQFT-type gluing
    properties. In this survey, we explain the formal structure and
    construction of bordered Floer homology and sketch how it can be
    used to compute some aspects of Heegaard Floer theory.
\end{abstract}


\maketitle 

\section{Introduction}

Heegaard Floer homology, introduced in a series of papers
\cite{OS04:HolomorphicDisks, OS04:HolDiskProperties, OS06:HolDiskFour}
of  Zolt\'an Szab\'o and the second author, has become a useful tool in
$3$- and $4$-dimensional topology. The Heegaard Floer invariants
contain subtle topological information, allowing one to detect the
genera of knots and homology classes \cite{OS04:GenusBounds};
detect fiberedness for knots~\cite{OS05:Contact,
  Ghiggini08:FiberedGenusOne, Ni07:FiberedKnot, Juhasz06:Sutured, Juhasz08:SuturedDecomp} and
$3$-manifolds~\cite{OS04:symplectic, Ni09:FiberedMfld,
  AiPeters10:twistedtorus, AiNi09:fiberedtorus}; bound the slice
genus~\cite{OS03:4BallGenus} and unknotting
number~\cite{OS05:unknotting, Owens08:unknotting}; prove tightness and obstruct Stein
fillability of contact structures~\cite{OS05:Contact,LiscaStipsicz09:contactseifert}; and more. It
has been useful for resolving a number of conjectures, particularly
related to questions about Dehn surgery \cite{OS05:surgeries,
  Greene:lensrealization}; see also~\cite{KMOS}. It is either known or conjectured to be
equivalent to several other gauge-theoretic or holomorphic curve
invariants in low-dimensional topology, including monopole Floer
homology~\cite{KronheimerMrowka}, embedded contact
homology~\cite{Hutchings02:ECH}, and the Lagrangian matching
invariants of 3- and
4-manifolds~\cite{Usher06:Vortices,Perutz:Matching1,Perutz:Matching2}. Heegaard
Floer homology is known to relate to Khovanov
homology~\cite{BrDCov,GrigsbyWehrli:detects,Hedden09:unknot},
and more relations with Khovanov-Rozansky type homologies are
conjectured~\cite{DGR06:Super}.

Heegaard Floer homology has several variants; the technically simplest
is $\HFa$, which is sufficient for most of the $3$-dimensional
applications discussed above. Bordered Heegaard Floer homology, the focus of
this paper, is an extension of $\HFa$ to $3$-manifolds with
boundary~\cite{LOT1}. This extension gives a conceptually satisfying
way to compute essentially all aspects of the Heegaard Floer
package related to $\HFa$. (There are also other algorithms
for computing many parts of Heegaard Floer theory
\cite{SarkarWang07:ComputingHFhat,LMW06:CombinatorialCobordismMaps,OSS:nice,MOS06:CombinatorialDescrip,MOST07:CombinatorialLink,MOT:grid,ManolescuOzsvath:surgery,OSS09:singular,OS09:cube}.)

We will start with the formal structure of bordered Heegaard
Floer homology. Most of the paper is then devoted to sketching its
definition. We conclude by explaining how bordered Floer homology can
be used for calculations of Heegaard Floer invariants.

\section{Formal structure}\label{sec:formal-structure}
\subsection{Review of Heegaard Floer theory}
Heegaard Floer theory has many components. Most basic among them, it associates:
\begin{itemize}
\item To a closed, connected, oriented $3$-manifold $Y$, an abelian
  group $\HFa(Y)$ and $\ZZ[U]$-modules $\HF^+(Y)$, $\HF^-(Y)$ and
  $\HF^\infty(Y)$. These
  are the homologies of chain complexes $\CFa(Y)$, $\CF^+(Y)$,
  $\CF^-(Y)$ and $\CF^\infty(Y)$, respectively. The chain complexes
  (and their homology groups) decompose into
  $\SpinC$-structures,
  $\CF(Y)=\bigoplus_{\spinc\in\SpinC(Y)}\CF(Y,\spinc)$, where $\CF$
  is any of the four chain complexes. Each $\CF(Y,\spinc)$ has a
  relative grading
  modulo the divisibility of $c_1(\spinc)$~\cite{OS04:HolomorphicDisks}.
  The chain complex $\CFa(Y)$ is the $U=0$ specialization of $\CF^-(Y)$.
\item To a smooth, compact, oriented cobordism $W$ from $Y_1$
  to~$Y_2$, maps $F_W\co \HF(Y_1)\to \HF(Y_2)$ induced by chain maps
  $f_W\co \CF(Y_1)\to \CF(Y_2)$.%
  \footnote{For $\CF^-$ and $\CF^\infty$, we mean the
    completions with respect to the formal variable~$U$.}
  These maps decompose according to $\SpinC$-structures
  on~$W$~\cite{OS06:HolDiskFour}.
\end{itemize}
\noindent
The maps $F_W$ satisfy a TQFT composition law: 
\begin{itemize}
\item 
  If $W'$ is another cobordism, from $Y_2$ to $Y_3$,
  then
  $
  F_{W'}\circ F_W=F_{W'\circ W}
  $~\cite{OS06:HolDiskFour}.
\end{itemize}

The Heegaard Floer invariants are defined by counting
pseudoholomorphic curves in symmetric products of Heegaard surfaces.
The Heegaard Floer groups were conjectured to be
equivalent to the monopole Floer homology groups (defined by counting solutions of the
Seiberg-Witten equations), via the correspondence:
\begin{align*}
  \HF^+(Y)&\longleftrightarrow\HMto(Y)\\
  \HF^-(Y)&\longleftrightarrow \HMfrom(Y)\\
  \HF^\infty(Y)&\longleftrightarrow\HMbar(Y),
\end{align*}
and similarly for the corresponding cobordism maps. A proof of this
conjecture has recently been announced by \c{C}a\u{g}atay
Kutluhan, Yi-Jen Lee and Clifford Taubes
\cite{KutluhanLeeTaubes:HFHMI,
  KutluhanLeeTaubes:HFHMII,KutluhanLeeTaubes:HFHMIII,KutluhanLeeTaubes:HFHMIV,KutluhanLeeTaubes:HFHMV}. Vincent Colin, Paulo
Ghiggini and Ko Honda have announced an independent proof for the $U=0$
specialization~\cite{ColinGhigginiHonda11:summary}.

In particular, the Heegaard Floer package contains enough information
to detect exotic smooth structures on
$4$-manifolds~\cite{OS04:symplectic, Mark:knottedsurf}. For closed
$4$-manifolds, this information is contained in $\HF^+$ and $\HF^-$;
the weaker invariant $\HFa$ is not useful for distinguishing smooth
structures on closed $4$-manifolds.

\subsection{The structure of bordered Floer theory}
Bordered Floer homology is an extension of $\HFa$ to $3$-manifolds
with boundary, in a TQFT form.  Bordered Floer homology associates:
\begin{itemize}
\item To a closed, oriented, connected surface $F$, together with some extra
  markings (see Definition~\ref{def:PMC}), a differential graded
  (\emph{dg}) algebra~$\Alg(F)$.
\item To a compact, oriented $3$-manifold $Y$ with connected boundary,
  together with a diffeomorphism $\phi\co F\to \bdy Y$ marking the
  boundary, a module over $\Alg(F)$. Actually, there are two
  different invariants for $Y$: $\CFDa(Y)$, a left \dg module over
  $\Alg(-F)$, and $\CFAa(Y)$, a right $\Ainf$ module over $\Alg(F)$,
  each well-defined up to quasi-isomorphism.
  We will sometimes refer to a $3$-manifold $Y$ with $\bdy Y=F$; we
  actually mean $Y$ together with an identification $\phi$ of $\bdy Y$
  with $F$. We call this data a \emph{bordered $3$-manifold.}
\item More generally, to a $3$-manifold $Y$ with two boundary
  components $\bdy_LY$ and $\bdy_RY$, diffeomorphisms $\phi_L\co F_L\to
  \bdy_LY$ and $\phi_R\co F_R\to \bdy_RY$ and a framed arc $\gamma$
  from $\bdy_LY$ to $\bdy_RY$ (compatible with $\phi_L$ and $\phi_R$
  in a suitable sense), a \dg
  bimodule $\CFDDa(Y)$ with commuting left actions of $\Alg(-F_L)$ and
  $\Alg(-F_R)$; an $\Ainf$ bimodule $\CFDAa(Y)$ with a left action of
  $\Alg(-F_L)$ and a right action of $\Alg(F_R)$; and an
  $\Ainf$ bimodule $\CFAAa(Y)$ with commuting right actions of $\Alg(F_L)$
  and $\Alg(F_R)$. Each of $\CFDDa(Y)$, $\CFDAa(Y)$ and $\CFAAa(Y)$ is
  well-defined up to quasi-isomorphism.
 We call the data $(Y,\phi_L,\phi_R,\gamma)$ a \emph{strongly
    bordered $3$-manifold with two boundary components}.
\end{itemize}
To keep the sidedness straight, note that type $D$ boundaries are
always on the left, and type $A$ boundaries are always on the
right; and
for $[0,1]\times F$,
the boundary component on the left side, $\{0\}\times F$,
is oriented as $-F$, while the one on the right side is oriented as~$F$.

Gluing $3$-manifolds corresponds to tensoring invariants; for any
valid tensor product (necessarily matching $D$ sides with $A$ sides)
there is a corresponding gluing. More
concretely~\cite{LOT1,LOT2}
\begin{itemize}
\item Given $3$-manifolds $Y_1$ and $Y_2$ with $\bdy Y_1=F=-\bdy Y_2$,
  there is a quasi-isomorphism 
  \begin{equation}
    \CFa(Y_1\cup_FY_2)\simeq
    \CFAa(Y_1)\DTP_{\Alg(F)}\CFDa(Y_2),\label{eq:pairing-1}  
  \end{equation}
where $\DTP$ denotes the
  derived tensor product. So, $\HFa(Y_1\cup_FY_2)\cong
  \Tor_{\Alg(F)}(\CFAa(Y_1),\CFDa(Y_2))$.
\item More generally, given $3$-manifolds $Y_1$ and $Y_2$ with $\bdy
  Y_1=-F_1\amalg F_2$ and $\bdy Y_2=-F_2\amalg F_3$ there are
  quasi-isomorphisms of bimodules corresponding to any valid tensor
  product. For instance:
  \begin{align*}
    \CFDDa(Y_1\cup_{F_2}Y_2)
     &\simeq\CFDAa(Y_1)\DTP_{\Alg(F_2)}\CFDDa(Y_2)\\
     &\simeq \CFDAa(Y_2)\DTP_{\Alg(-F_2)}\CFDDa(Y_1).
  \end{align*}
  The mnemonic is that any tensor product which matches
  $A$'s with~$D$'s corresponds to a valid gluing.
  Also, $F_1$ or $F_3$ may be $S^2$ (or empty), in which case these
  statements reduce to pairing theorems for a module and a
  bimodule.
  If both $F_1$ and $F_3$ are empty, these reduce to
  Equation~[\ref{eq:pairing-1}].
\end{itemize}
We refer to
theorems of this kind as \emph{pairing theorems}. 

\begin{remark}\label{remark:efficient}
  One feature of the theory is that if one uses a suitable model of
  the derived tensor product then these pairing theorems are
  efficient, in the sense that $\CFa(Y_1\cup_F Y_2)$ and
  $\CFAa(Y_1)\DTP_{\Alg(F)}\CFDa(Y_2)$ have roughly the same number of
  generators.
\end{remark}

There is also a self-pairing theorem. Let $Y$ be a $3$-manifold with
$\bdy Y=-F\amalg F$ and $\gamma$ be a framed arc connecting
corresponding points in the boundary components of~$Y$.
Let $Y^\circ$
denote the \emph{generalized open book associated to $(Y,\gamma)$}
obtained by gluing the two boundary components of $Y$ and performing
framed surgery along the $S^1$ in the result corresponding to
$\gamma$. Let $\gamma^\circ$ be the null-homologous knot in $Y^\circ$
corresponding to $\gamma$. Then:
\begin{itemize}
\item 
  The Hochschild homology of $\CFDAa(Y)$ is given~\cite{LOT2} by
\[
\HH_*(\CFDAa(Y))\cong \HFKa(Y^\circ,\gamma^\circ).
\]
\end{itemize}

These invariants satisfy a number of duality
properties~\cite{LOTHomPair}, e.g.:
\begin{itemize}
\item The algebra $\Alg(F)$ is the opposite algebra of $\Alg(-F)$.
  (There are also more subtle duality properties of the algebras;
  see Remark~\ref{remark:koszul-duality}.)
\item The module $\CFDa(Y)$ is dual (over $\Alg(F)$) to $\CFAa(-Y)$:
  \begin{align*}
  \CFDa(Y)&\simeq \Mor_{\Alg(-F)}(\CFAa(-Y),\Alg(-F))\\
  \CFAa(Y)&\simeq \Mor_{\Alg(F)}(\CFDa(-Y),\Alg(F)).
  \end{align*}
  (Here, we are
  using the fact that $\Alg(-F)=\Alg(F)^\op$ to exchange left actions
  by $\Alg(-F)$ and right actions by $\Alg(F)$.)
\item The module $\CFDDa(-Y)$ is the one-sided dual of $\CFAAa(Y)$:
  \[
  \Mor_{\Alg(-F_2)}(\CFDDa(-Y),\Alg(-F_2))\simeq \CFAAa(Y). 
  \]
  The symmetric statement also holds, as does the corresponding statement for $\CFDAa(Y)$.
\item Given a strongly bordered $3$-manifold $Y$ with two boundary
  components, let $\tau_\bdy^{-1}(Y)$ denote the same manifold $Y$ but with the framing on the framed arc $\gamma$ increased by $1$. Then the two-sided dual of $\CFDDa(Y)$ is $\CFAAa(-\tau_\bdy(Y))$, i.e., 
  \begin{multline*}
  \Mor_{\Alg(-F_1)\otimes\Alg(-F_2)}(\CFDDa(Y),\Alg(-F_1)\otimes\Alg(-F_2))\\
    \simeq \CFAAa(-\tau_\bdy^{-1}(Y)).
  \end{multline*}
\end{itemize}
(We include the last property here
mainly to indicate that some caution is needed.)

As a consequence of these dualities, one can give pairing theorems
using the $\Hom$ functor rather than the tensor
product~\cite{LOTHomPair}, e.g.:
\begin{itemize}
\item Let $Y_1$ and $Y_2$ be $3$-manifolds with $\bdy Y_1=\bdy Y_2=F$. Then 
  \begin{align*}
  \CFa(-Y_1\cup_\bdy Y_2)
   &\simeq \Mor_{\Alg(-F)}(\CFDa(Y_1),\CFDa(Y_2))\\
   &\simeq \Mor_{\Alg(F)}(\CFAa(Y_1),\CFAa(Y_2)).
  \end{align*}
  Similar statements hold 
  for bimodules;
  for instance, if $Y_1$ has another boundary component $F_0$ and $Y_2$ has another boundary component $F_2$ then
  \begin{align*}
  \CFAAa(-Y_1\cup_\bdy Y_2)
    &\simeq \Mor_{\Alg(-F)}(\CFDDa(Y_1),\CFDAa(Y_2))\\
    &\simeq \Mor_{\Alg(F)}(\CFDAa(Y_1),\CFAAa(Y_2)).
  \end{align*}
\item Given $3$-manifolds $Y_1$ and $Y_2$ with $\bdy Y_1=F=-\bdy Y_2$, 
  \begin{multline}\label{eq:MorPair1}
  \CFa(Y_1\cup_FY_2) \simeq
  \Mor_{\Alg(-F)\otimes\Alg(F)}(\CFDDa([0,1]\times F),\\
    \CFDa(Y_1)\otimes_\Field\CFDa(Y_2)).
  \end{multline}
Similarly, if $Y_2$ has another boundary component $F'$ then 
  \begin{multline}\label{eq:MorPair2}
  \CFDa(Y_1\cup_FY_2)\simeq
  \Mor_{\Alg(-F)\otimes\Alg(F)}(\CFDDa([0,1]\times F),\\
    \CFDa(Y_1)\otimes_\Field\CFDDa(Y_2)).
  \end{multline}
  (If both $Y_1$ and $Y_2$ had
  two boundary components then the left hand side would pick up a
  change of framing.)
\end{itemize}

\begin{remark}
  Some of the duality properties above can also be seen from
  the Fukaya\hyp categorical perspective~\cite{AurouxBordered}.
\end{remark}

\begin{remark}
  The pairing theorems using homomorphisms are not usually efficient
  in the sense of Remark~\ref{remark:efficient}.
\end{remark}

\begin{remark}
  It is conceptually cleaner to think only about type \DA\
  bimodules. However, as their definitions involve simpler moduli
  spaces of curves, it is usually easier to compute type \DD\
  bimodules.
\end{remark}
\begin{remark}
  It is natural to expect that to a $4$-manifold with corners one
  would associate a map of bimodules, satisfying certain gluing
  axioms. We have not done this; however, as discussed below, even
  without this bordered
  Floer homology allows one to compute the maps $\hat{F}_W$ associated to
  cobordisms $W$ between closed $3$-manifolds.
\end{remark}

\section{The algebras}\label{sec:algebras}
As mentioned earlier, the bordered Floer algebras are associated to
surfaces together with some extra markings. We encode these markings
as \emph{pointed matched circles} $\PMC$, which we discuss next. We then
introduce a simpler algebra, $\Alg(n)$, depending only on an integer
$n$, of which the bordered Floer algebras $\Alg(\PMC)$ are
subalgebras. The definition of $\Alg(\PMC)$ itself is given in the
last subsection.
\subsection{Pointed matched circles}\label{sec:pmc}
\begin{definition}\label{def:PMC}
  A \emph{pointed matched circle} $\PMC$ consists of:
  \begin{itemize}
  \item an oriented circle $Z$,
  \item $4k$ points $\CircPts=\{a_1,\dots,a_{4k}\}$ in $Z$,
  \item a matching $M$ of the points in $\CircPts$ in pairs, which we
    view as a fixed-point free involution $M\co \CircPts\to \CircPts$,
    and
  \item a basepoint $z\in Z\setminus\CircPts$.
  \end{itemize}
  We require that performing surgery on $Z$ along the matched pairs of
  points yields a connected $1$-manifold.
\end{definition}

A pointed matched circle $\PMC$ with $|\CircPts|=4k$ specifies:
\begin{itemize}
\item A closed surface $F(\PMC)$ of genus $k$, as follows. Fill $Z$
  with a disk~$D$. Attach a $2$-dimensional $1$-handle to each pair of
  points in
  $\CircPts$ matched by~$M$. By hypothesis, the result has connected
  boundary; fill that boundary with a second disk.
\item A distinguished disk in $F(\PMC)$: the disk $D$ (say).
\item A basepoint $z$ in the boundary of the distinguished disk.
\end{itemize}

\begin{figure}
  \centering
  \includegraphics[scale=.6666667]{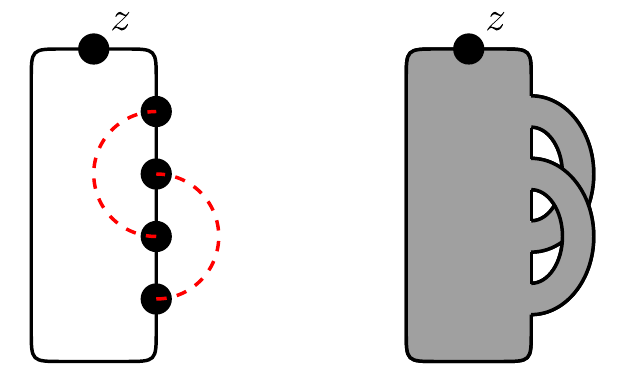} 
  \caption{\textbf{A pointed matched circle (left) and the surface it specifies (right).} On the left, the dashed lines indicate the matching. On the right, a copy of $\bD^2$ is glued to the boundary of the surface shown.}
  \label{fig:pmc-to-surf}
\end{figure}

See Figure~\ref{fig:pmc-to-surf} for an example.

\begin{remark}
  Different pointed matched circles may specify homeomorphic
  surfaces. If $F(\PMC)\cong F(\PMC')$ then, while the bordered Floer algebras
  $\Alg(\PMC)$ and $\Alg(\PMC')$ will typically be different (cf.~Example~\ref{eg:algebras}), their
  module categories will be equivalent; this follows from formal
  properties of the bordered Floer package (pairing theorems and the
  fact that the bimodule $\CFDAa([0,1]\times F(\PMC))$ is homotopy
  equivalent to $\Alg(\PMC)$).
\end{remark}

\begin{remark}
  Matched circles can be seen as a special case of fat
    graphs~\cite{Penner87:MCGroupoid}. They are also dual to the
  typical representation of a genus~$g$ surface
  as a $4g$-gon with sides glued together.
\end{remark}

\subsection{The strands algebra}\label{sec:strands-algebra}
We next define a differential algebra $\Alg(n)$, depending only on an
integer $n$; the algebra $\Alg(\PMC)$ associated to a pointed matched
circle with $|\CircPts|=4k$ will be a subalgebra of $\Alg(4k)$.

The algebra $\Alg(n)$ has an $\Field$-basis consisting of all triples
$(S,T,\phi)$ where $S$ and $T$ are subsets of
$\underline{n}=\{1,\dots,n\}$ and $\phi\co S\to T$ is a bijection such
that for all $s\in S$, $\phi(s)\geq s$. Given such a map~$\phi$, let
$\Inv(\phi)=\{(s_1,s_2)\in S\times S\mid s_1<s_2,\
\phi(s_2)<\phi(s_1)\}$ and $\inv(\phi)=|\Inv(\phi)|$, so $\inv(\phi)$
is the number of \emph{inversions} of~$\phi$.

The product $(S,T,\phi)\cdot (U,V,\psi)$ in $\Alg(n)$ is defined to be
$0$ if $U\neq T$ or if $U=T$ but $\inv(\psi\circ \phi)\neq
\inv(\psi)+\inv(\phi)$. If $U=T$ and $\inv(\psi\circ \phi)=
\inv(\psi)+\inv(\phi)$ then let $(S,T,\phi)\cdot
(U,V,\psi)=(S,V,\psi\circ\phi)$.  In particular, the elements
$(S,S,\Id)$ (where $\Id$ denotes the identity map) are the
indecomposable idempotents in $\Alg(n)$.

Given a generator $(S,T,\phi)\in\Alg(n)$ and an element
$\sigma=(s_1,s_2)\in\Inv(\phi)$, let $\phi_\sigma\co S\to T$ be the
map defined by $\phi_\sigma(s)=\phi(s)$ if $s\neq s_1,s_2$;
$\phi_\sigma(s_1)=\phi(s_2)$; and $\phi_\sigma(s_2)=\phi(s_1)$. Define
a differential on $\Alg(n)$ by:
\[
\bdy(S,T,\phi)=\!\!\sum_{
  \begin{subarray}{c}
    \sigma\in\Inv(\phi)\\
    \inv(\phi_\sigma)=\inv(\phi)-1
  \end{subarray}}\!\!
  (S,T,\phi_\sigma).
\]

\begin{figure*}
  \centering
  \input{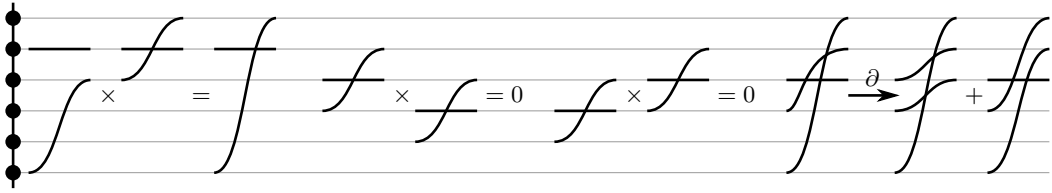}
  \caption{\textbf{The strands algebra.} A product, two vanishing
    products, and a differential. In this notation, the restrictions
    on the number of inversions means that elements with double
    crossings in the product or differential are set to $0$.}
  \label{fig:strands-alg}
\end{figure*}

See Figure~\ref{fig:strands-alg} for a graphical representation.

Given a generator $(S,T,\phi)\in\Alg(n)$, define the \emph{weight of
  $(S,T,\phi)$} to be the cardinality of $S$. 
Let $\Alg(n,i)$ be
the subalgebra of $\Alg(n)$ generated by elements of weight $i$, so
$\Alg(n)=\bigoplus_{i=0}^n\Alg(n,i)$.

\subsection{The algebra associated to a pointed matched
  circle}\label{sec:alg-of-pmc}
Fix a pointed matched circle $\PMC=(Z,\CircPts,M,z)$ with
$|\CircPts|=4k$. After cutting $Z$ at~$z$, the orientation of $Z$ identifies
$\CircPts$ with $\underline{4k}$, so we can view $M$ as a matching of~$\underline{4k}$.

Call a basis element $(S,T,\phi)$ of $\Alg(4k)$ \emph{equitable} if
no two elements of $\underline{4k}$ that are matched (with respect to $M$)
both occur in~$S$, and no two elements of $\underline{4k}$ that are matched
both occur in~$T$.

Given equitable basis elements $x=(S,T,\phi)$ and
$y=(S',T',\allowbreak\psi)$ of $\Alg(4k)$, we say that $x$ and $y$ are related by
\emph{horizontal strand swapping}, and write $x\sim y$, if there is a
subset $U\subset S$ such that $S'=(S\setminus U)\cup M(U)$,
$\phi|_{S\setminus U}=\psi|_{S\setminus U}$,
$\phi|_{U}=\Id_U$ and $\psi|_{M(U)}=\Id_{M(U)}$.

\begin{figure}
  \centering
  \input{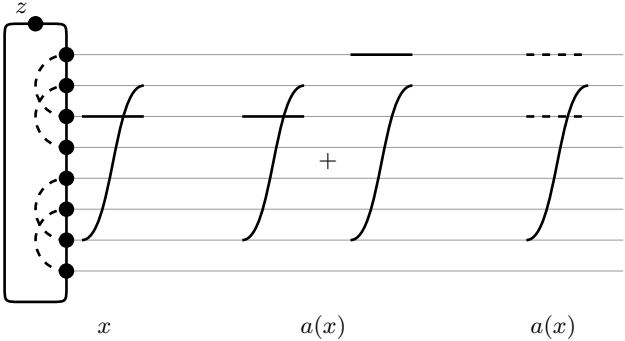}
  \caption{\textbf{The algebra $\Alg(\PMC)$.} An example of the operation $a(x)$ and a shorthand for the resulting element of $\Alg(\PMC)\subset \Alg(4k)$.}
  \label{fig:Alg-of-pmc}
\end{figure}

Given an equitable basis element $x$ of $\Alg(4k)$, let
$a(x)=\sum_{y\sim x}y$.
See Figure~\ref{fig:Alg-of-pmc} for an example.
Define $\Alg(\PMC)\subset\Alg(4k)$ to be the $\Field$-subspace with
basis $\{a(x)\mid x\text{ is equitable }\}$. It is straightforward to
verify that $\Alg(\PMC)$ is a differential subalgebra of
$\Alg(4k)$. We will call the elements $a(x)$ \emph{basic generators of $\Alg(\PMC)$}.
If $x$ is not equitable, set $a(x)=0$, and extend
$a$ linearly to a map $a\co \Alg(4k)\to \Alg(\PMC)$. (This is not
an algebra homomorphism.)

Indecomposable idempotents of $\Alg(\PMC)$ correspond to subsets of
the set of matched pairs in~$\CircPts$.  These generate a subalgebra
$\Idem(\PMC)$ where all strands are horizontal.
The algebra $\Alg(\PMC)$
decomposes as $\bigoplus_{i=-k}^k\Alg(\PMC,i)$ where $\Alg(\PMC,i) =
\Alg(\PMC) \cap \Alg(4k,k+i)$.

As the figures suggest, we often think of elements of $\Alg(\PMC)$ in
terms of sets of {\em chords} in $(Z,\CircPts)$, i.e., arcs in $Z$ with
endpoints in $\CircPts$, with orientations induced from $Z$. Given a
chord $\rho$ in
$(Z\setminus\{z\},\CircPts)$ let $\rho^-$ (respectively $\rho^+$) be
the initial (respectively terminal) endpoint of $\rho$.  Given a set
$\rhos=\{\rho_i\}$ of chords in $(Z\setminus \{z\},\CircPts)$ such
that no two $\rho_i\in\rhos$ have the same initial (respectively
terminal) endpoint, let $\rhos^-=\{\rho_i^-\}$ and
$\rhos^+=\{\rho_i^+\}$; we can think of $\rhos$ as a map
$\phi_\rhos\co \rhos^-\to \rhos^+$. 
Let
\[
  a(\rhos)=\!\!\!\!\sum_{
    S_0 \cap \rhos^- = S_0 \cap \rhos^+ = \emptyset
  }\!\!\!\!
  (S_0 \cup \rhos^-,S_0 \cup
  \rhos^+,\phi_\rhos\amalg\Id_{S\setminus\rhos^-}).
\]

\begin{example}\label{eg:algebras}
  The algebra associated to the unique pointed matched circle for
  $S^2$ is $\Field$. The algebra $\Alg(T^2,0)$ associated to the
  unique pointed matched circle for $T^2$, 
  with $(1 \leftrightarrow 3, 2 \leftrightarrow 4)$,
  is given by
  the path algebra with relations:
  \[
  \mathcenter{
    \begin{tikzpicture}
      \node at (0,0) (I1) {$I_1$};
      \node at (2,0) (I2) {$I_2$};
      \draw[->, bend left=30] (I1) to node[above]{$\rho_{1,2},\ \rho_{3,4}$} (I2);
      \draw[->, bend left=30] (I2) to node[below]{$\rho_{2,3}$} (I1);
    \end{tikzpicture}}
  \mathcenter{\Big/(\rho_{2,3}\rho_{1,2}=\rho_{3,4}\rho_{2,3}=0).}
  \]
  The algebra associated to the pointed matched circle $\PMC$
  for a genus $2$ surface with matching 
  $
  (1\leftrightarrow3,\ 
  2\leftrightarrow4,\ 
  5\leftrightarrow7,\ 
  6\leftrightarrow8)
  $
  has Poincar\'e polynomial
  \cite[Section 4]{LOT2}
  \[
  \sum_i \dim_\Field H_*(\Alg(\PMC,i))T^i=T^{-2}+32T^{-1}+98+32T^1+T^2.
  \] 
  The algebra associated to the pointed
  matched circle $\PMC'$ for a genus $2$ surface with matching
  $
  (1\leftrightarrow 5,\ 
  2\leftrightarrow 6,\ 
  3\leftrightarrow 7,\ 
  4\leftrightarrow 8)
  $
  has Poincar\'e polynomial
  \[
  \sum_i \dim_\Field H_*(\Alg(\PMC',i))T^i=T^{-2}+32T^{-1}+70+32T^1+T^2.
  \]
\end{example}

  The ranks in the genus two examples which are equal are explained by
  the observations that for any pointed matched circle,
  $\Alg(\PMC,-k)\cong \Field$; $\Alg(\PMC,-k+1)$ has no differential;
  the dimension of $\Alg(\PMC,-k+1)$ is independent of the matching; and 
  the following.

\begin{remark}\label{remark:koszul-duality} The algebras
  $\Alg(\PMC,i)$ and $\Alg(-\PMC,-i)$ are Ko\-szul dual. (Here, $-\PMC$
  denotes the pointed matched circle obtained by reversing the
  orientation on $Z$.) Also, given a pointed matched circle $\PMC$ for
  $F$, let $\PMC_*$ denote the pointed matched circle corresponding to
  the dual handle decomposition of $F$. Then $\Alg(\PMC,i)$ and
  $\Alg(\PMC_*,i)$ are Koszul dual. In particular, $\Alg(\PMC,-i)$ is
  quasi-isomorphic to $\Alg(\PMC_*,i)$~\cite{LOTHomPair}.
\end{remark}

\begin{remark}
  In Rumen Zarev's bordered-sutured extension of the
  theory~\cite{Zarev09:BorSut}, the strands algebra $\Alg(n,k)$ has a
  topological interpretation as the algebra associated to a disk with
  boundary
  sutures.
\end{remark}

\section{Combinatorial representations of bordered
  3-manifolds}\label{sec:heegaard-diagrams}
A \emph{bordered $3$-manifold} is a $3$-manifold $Y$ together with an
orientation-preserving homeomorphism $\phi\co F(\PMC)\to \bdy Y$ for
some pointed matched circle $\PMC$. Two bordered $3$-manifolds
$(Y_1,\,\phi_1\co F(\PMC_1)\to \bdy Y_1)$ and $(Y_2,\,\phi_2\co
F(\PMC_2)\to \bdy Y_2)$ are called \emph{equivalent} if there is an
orientation-preserving
homeomorphism $\psi \co Y_1\to Y_2$ such that $\phi_1=\phi_2\circ
\psi$; in particular, this implies that $\PMC_1=\PMC_2$. Bordered
Floer theory associates homotopy equivalence classes of modules to
equivalence classes of bordered $3$-manifolds.
Just as the bordered Floer algebras are associated to combinatorial
representations of
surfaces, not directly to
surfaces, the bordered Floer modules are associated to
combinatorial representations of bordered $3$-manifolds.

\subsection{The closed case}
Recall that a \emph{$3$-dimensional handlebody} is a regular neighborhood of a
connected graph in $\RR^3$.
These are a particularly simple class of
$3$-manifolds with boundary.  
According to a classical result of
Poul Heegaard~\cite{Heegaard98}, every closed, orientable $3$-manifold
can be obtained as a union of two such
handlebodies, $H_\alpha$ and $H_\beta$. Such a representation is
called a {\em Heegaard splitting}. A Heegaard splitting along
an orientable surface~$\Sigma$ of genus~$g$ can be
represented by a {\em Heegaard diagram}: a pair of $g$-tuples of
pairwise disjoint, homologically linearly independent, 
embedded circles 
$\alphas=\{\alpha_1,\dots,\alpha_g\}$ and
$\betas=\{\beta_1,\dots,\beta_g\}$ in~$\Sigma$.  These curves are
chosen so that each $\alpha_i$ (resp. $\beta_i$) bounds a disk in
the handlebody $H_\alpha$ (resp. $H_\beta$).
Any two Heegaard diagrams for the same manifold~$Y$ are related by a
sequence of moves, called \emph{Heegaard
  moves}; see, for instance,~\cite{Scharlemann02:HeegaardSurvey} or~\cite[Section~2.1]{OS04:HolomorphicDisks}.

\subsection{Representing 3-manifolds with boundary}
The story extends easily to $3$-manifolds with boundary, using
a slight generalization of handlebodies. 
A \emph{compression body} (with both boundaries connected) is the result of
starting with a connected orientable surface~$\Sigma_2$ times $[0,1]$ and
then attaching thickened disks ($3$-dimensional $2$-handles) along
some number of homologically linearly independent, disjoint circles in
$\Sigma_2\times\{0\}$. A
compression body has two boundary components, $\Sigma_1$ and
$\Sigma_2$, with genera $g_1\leq g_2$. Up to homeomorphism, a
compression body is determined by its boundary.

A \emph{Heegaard decomposition} of a $3$-manifold $Y$ with non\-empty, connected boundary
is a decomposition $Y=H_\alpha\cup_{\Sigma_2} H_\beta$ where $H_\alpha$ is a
compression body and $H_\beta$ is a handlebody. Let $g$ be the genus of
$\Sigma_2$ and $k$ the genus of $\bdy Y$. A \emph{Heegaard diagram}
for $Y$ is gotten by choosing $g$ pairwise disjoint circles
$\beta_1,\dots,\beta_g$ in $\Sigma_2$ and $g-k$ disjoint circles
$\alpha^c_1,\dots,\alpha^c_{g-k}$ in $\Sigma_2$ so that:
\begin{itemize}
\item The circles $\beta_1,\dots,\beta_g$ bound disks
  $D_{\beta_1},\dots,D_{\beta_g}$ in $H_\beta$ such that
  $H_\beta\setminus(D_{\beta_1}\cup\dots\cup D_{\beta_g})$ is
  topologically a ball, and
\item The circles $\alpha_1^c,\dots,\alpha_{g-k}^c$ bound disks
  $D_{\alpha_1^c},\allowbreak\dots,\allowbreak D_{\alpha_{g-k}^c}$ in $H_\alpha$ such that
  $H_\alpha\setminus (D_{\alpha_1^c}\cup\dots\cup D_{\alpha_{g-k}^c})$ is
  topologically the product of a surface and an interval.
\end{itemize}

\begin{figure}
  \centering
  \input{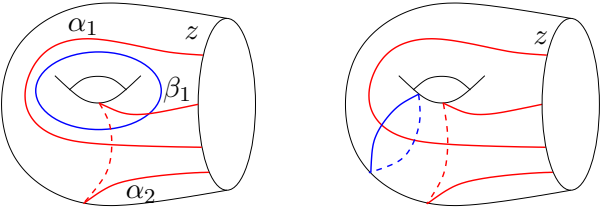}
  \caption{\textbf{Framed handlebodies.}
    Two different framings of a genus one handlebody are shown. In
    these examples, $g=k=1$; in particular, there are no $\alpha$-circles.}
  \label{fig:Handlebodies}
\end{figure}

To specify a parametrization, or bordering, of~$\partial Y$, we
need a little more data. A \emph{bordered
  Heegaard diagram} for $Y$ is a tuple 
\[
\HD=(\overline{\Sigma},\overbrace{\alpha_1^c,\dots,\alpha_{g-k}^c}^{\alphas^c},\overbrace{\alpha_1^a,\dots,\alpha_{2k}^a}^{\alphas^a},\overbrace{\beta_1,\dots,\beta_g}^\betas,
z)
\]
where:
\begin{itemize}
\item $\overline{\Sigma}$ is an oriented surface with a single
  boundary component.
\item
  $(\overline{\Sigma}\cup_\bdy \bD^2,\alphas^c,\betas)$
  is a Heegaard diagram for~$Y$.
\item $\alpha_1^a,\dots,\alpha_{2k}^a$ are pairwise disjoint, embedded
  arcs in $\overline{\Sigma}$ with boundary on
  $\bdy\overline{\Sigma}$, and are disjoint from the $\alpha_i^c$.
\item  $\overline{\Sigma}\setminus
  (\alpha_1^c\cup\dots\cup\alpha_{g-k}^c\cup \alpha_1^a\cup
  \dots\cup\alpha_{2k}^a)$ is a disk with $2(g-k)$ holes.
\item $z$ is a point in $\bdy\overline{\Sigma}$, disjoint from all of
  the $\alpha_i^a$.
\end{itemize}
See Figure~\ref{fig:Handlebodies} for two examples of bordered
Heegaard diagrams for the solid torus.

Let $\alphas=\alphas^a\cup\alphas^c$.

A bordered Heegaard diagram $\HD$ specifies a pointed matched circle
\[
\PMC(\HD)=(Z=\bdy\overline{\Sigma},\CircPts=(\alphas^a\cap\bdy\overline{\Sigma}),
M, z),
\]
where two points in $\CircPts$ are matched in $M$ if they lie on the
same $\alpha_i^a$. A bordered Heegaard diagram for $Y$ also specifies
an identification $\phi \co F(\PMC)\to \bdy Y$, well-defined up to
isotopy.

There are moves, analogous to Heegaard moves, relating any two
bordered Heegaard diagrams for equivalent bordered $3$-manifolds.

\section{The modules and bimodules}\label{sec:modules}
As discussed above, there are two
invariants of a $3$\hyp manifold $Y$ with boundary
$F(\PMC)$. $\CFDa(Y)$ has a straghtforward module structure but a
differential which counts holomorphic curves, while $\CFAa(Y)$
uses holomorphic curves to define the module structure itself. 

Fix a bordered Heegaard diagram $\HD=(\bSigma,\alphas,\betas,z)$ for
$Y$. Let $\Gen(\HD)$ be the set of $g$-tuples
$\x=\{x_i\}_{i=1}^g\subset\alphas\cap\betas$ so that there is exactly
one point $x_i$ on each $\beta$-circle and on each $\alpha$-circle and
there is at most one $x_i$ on each $\alpha$-arc. The invariant
$\CFAa(Y)$ is a direct sum of copies of $\Field$, one for each element
of $\Gen(\HD)$, while $\CFDa(Y)$ is a direct sum of elementary
projective $\Alg(-\bdy \HD)$-modules, one for each element of
$\Gen(\HD)$. Let $X(\HD)$ be the
$\Field$-vector space generated by $\Gen(\HD)$, which is also the
vector space underlying $\CFAa(Y)$.

Each generator $\x\in\Gen(\HD)$ 
determines a $\SpinC$-structure $\spinc(\x)\in\SpinC(Y)$;
the construction~\cite{LOT1} is an easy adaptation of the closed case
\cite[Section~2.6]{OS04:HolomorphicDisks}.

Before continuing to describe the bordered Floer modules, we digress
to briefly discuss the moduli spaces of holomorphic curves.

\subsection{Moduli spaces of holomorphic cuves}\label{sec:moduli-spaces}
Fix a bordered Heegaard diagram
$\HD=(\overline{\Sigma},\alphas,\betas,z)$. Let
$\Sigma=\overline{\Sigma}\setminus\bdy\overline{\Sigma}$. Choose a
symplectic form $\omega_\Sigma$ on $\Sigma$ giving it a cylindrical
end and a complex structure $j_\Sigma$ compatible with $\omega$,
making $\Sigma$ into a punctured Riemann surface. Let $p$ denote the
puncture in $\Sigma$.  We choose the $\alpha_i^a$ so that
their intersections with $\Sigma$ (also denoted $\alpha_i^a$)
are cylindrical ($\RR$-invariant) in a neighborhood of $p$.

We will consider curves
\[
u\co (S,\bdy S)\to
(\Sigma\times[0,1]\times\RR,\alphas\times\{1\}\times\RR\cup\betas\times\{0\}\times\RR)
\]
holomorphic with respect to an appropriate almost complex
structure~$J$, satisfying conditions spelled out in~\cite{LOT1}.
The reader may wish to simply think of a product complex structure
$j_\Sigma\times j_\bD$, though these complex structures may not be
  general enough to achieve transversality.

Such holomorphic curves $u$ have asymptotics in three places:
\begin{itemize}
\item $\Sigma\times[0,1]\times\{\pm\infty\}$. We consider curves
  asymptotic to $g$-tuples of strips $\x\times[0,1]\times\RR$ at
  $-\infty$ and $\y\times[0,1]\times\RR$ at $+\infty$, where
  $\x,\y\in\Gen(\HD)$.
\item $\{p\}\times[0,1]\times\RR$, which we denote $e\infty$. We
  consider curves asymptotic to chords $\rho_i$ in
  $(\bdy\bSigma,\CircPts)$ at a point
  $(1,t_i)\in[0,1]\times\RR$. (These are chords for the coisotropic
  foliation of $\bdy\overline{\Sigma}\times[0,1]\times\RR$, whose
  leaves are the circles $\bdy\overline{\Sigma}\times\{(s_0,t_0)\}$.)
  We impose the condition that these chords $\rho_i$ not cross $z\in
  \bdy\bSigma$.
\end{itemize}

Topological maps of this form can be grouped into homology classes. Let
$\pi_2(\x,\y)$ denote the set of homology classes of maps asymptotic
to $\x\times[0,1]$ at $-\infty$ and $\y\times[0,1]$ at $+\infty$. Then $\pi_2(\x,\x)$ is
canonically isomorphic to $H_2(Y,\bdy Y)$; $\pi_2(\x,\y)$ is nonempty
if and only if $\spinc(\x)=\spinc(\y)$; and if $\spinc(\x)=\spinc(\y)$
then $\pi_2(\x,\y)$ is an affine copy of $H_2(Y,\bdy Y)$, under
concatenation by elements of $\pi_2(\x,\x)$ (or
$\pi_2(\y,\y)$)~\cite{LOT1}. (Again, these results are easy adaptations of the
corresponding results in the closed case
\cite[Section~2]{OS04:HolomorphicDisks}.) Note that our usage of
$\pi_2(\x,\y)$ differs from the usage in~\cite{OS04:HolomorphicDisks},
where homology classes are allowed to cross $z$, but agrees with the usage in~\cite{LOT1}.

Given generators $\x,\y\in\Gen(\HD)$, a homology class
$B\in\pi_2(\x,\y)$ and a sequence
$\vec{\rhos}=(\rhos_1,\dots,\rhos_n)$ of sets
$\rhos_i=\{\rho_{i,1},\dots,\allowbreak\rho_{i,m_i}\}$ of Reeb chords, let
\[
\tcM^B(\x,\y;\vec{\rhos})
\]
denote the moduli space of embedded holomorphic curves $u$ in the
homology class $B$, asymptotic to $\x\times[0,1]$ at $-\infty$,
$\y\times[0,1]$ at $+\infty$, and $\rho_{i,j}\times(1,t_i)$ at
$e\infty$, for some sequence of heights $t_1<\cdots<t_n$. There is an
action of $\RR$ on $\tcM^B(\x,\y;\vec{\rhos})$, by translation. Let
$\cM^B(\x,\y;\vec{\rhos})=\tcM^B(\x,y;\vec{\rhos})/\RR$.

\begin{figure}
  \centering
  \includegraphics[scale=.57142857]{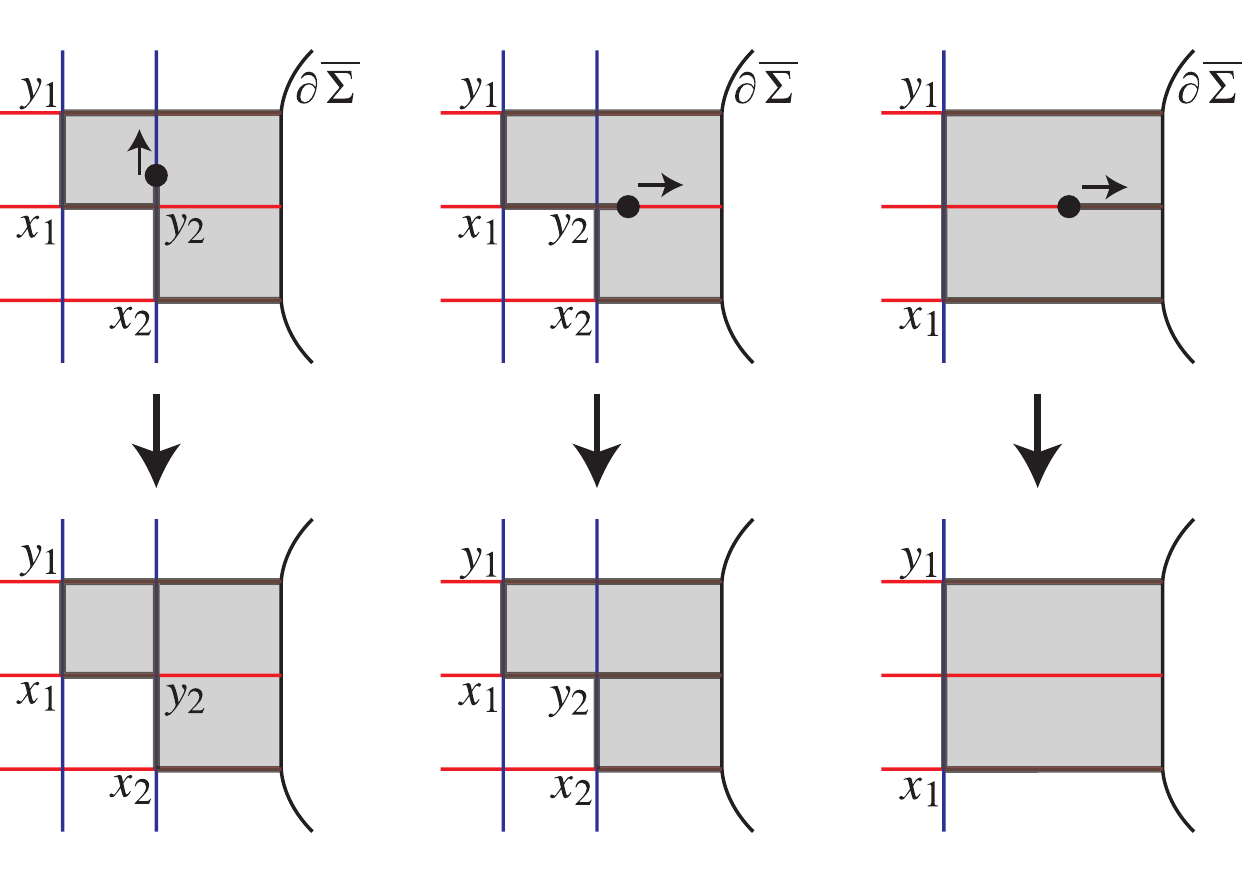} 
  \caption{\textbf{Degenerations of holomorphic curves.} Degenerations of types~\ref{item:two-story}, \ref{item:join} and~\ref{item:split} are shown, in that order. The dots indicate branch points, which can be thought of as the ends of cuts. (This figure is drawn from~\cite{LOT1}.)}
  \label{fig:degens-first-three}
\end{figure}

The modules $\CFDa(\HD)$ and $\CFAa(\HD)$ will be defined using counts
of $0$-dimensional moduli spaces $\cM^B(\x,\y;\vec{\rhos})$. Proving
that these modules satisfy $\bdy^2=0$ and the $\Ainf$ relations,
respectively, involves studying the ends of $1$-dimensional moduli
spaces. These ends correspond to the following four kinds of
degenerations:
\begin{enumerate}[label=D\arabic*.,ref=D\arabic*]
\item\label{item:two-story} Breaking into a two-story holomorphic
  building. That is, the $\RR$-coordinate of some parts of the curve
  go to $+\infty$ with respect to other parts, giving an element of
  $\cM^{B_1}(\x,\y;\vec{\rhos}_1)\times\cM^{B_2}(\x,\y;\vec{\rhos}_2)$
  where $B$ is the concatenation $B_1*B_2$ and $\vec\rhos$ is the
  concatenation $(\vec\rhos_1,\vec\rhos_2)$.
\item\label{item:join} Degenerations in which a boundary branch point
  of the projection $\pi_\Sigma\circ u$ approaches $e\infty$, in such
  a way that some chord $\rho_{i,j}$ splits into a pair of chords
  $\rho_a,\rho_b$ with $\rho_{i,j}=\rho_a\cup \rho_b$.  This degeneration results
  in a curve at $e\infty$, a \emph{join curve}, and an element of
  $\cM^B(\x,\y;\vec{\rhos}')$ where
  $\vec{\rhos}'$ is obtained by replacing the chord
  $\rho_{i,j}\in\rhos_i\in\vec{\rhos}$ with two chords, $\rho_a$ and~$\rho_b$.
\item\label{item:split} The difference in $\RR$-coordinates
  $t_{i+1}-t_i$ between two consecutive sets of chords $\rhos_i$ and
  $\rhos_{i+1}$ in $\vec{\rhos}$ going to $0$. In the process, some
  boundary branch points of $\pi_\Sigma\circ u$ may approach
  $e\infty$, degenerating a \emph{split curve}, along with an element of
  $\cM^B(\x,\y;\vec{\rhos}')$ where
  $\vec{\rhos}'=(\rhos_1,\dots,\rhos_i\uplus\rhos_{i+1},\dots,\rhos_n)$
  and $\rhos_i\uplus\rhos_{i+1}$ is gotten from
  $\rhos_i\cup\rhos_{i+1}$ by gluing together any pairs of chords
  $(\rho_{i,j},\rho_{i+1,\ell})$ where $\rho_{i,j}$ ends at the
  starting point of $\rho_{i+1,\ell}$ (i.e., $\rho_{i,j}^+=\rho_{i+1,\ell}^-$).
\item\label{item:shuffle} Degenerations in which a pair of boundary
  branch points of $\pi_\Sigma\circ u$ approach $e\infty$, causing a
  pair of chords $\rho_{i,j}$ and $\rho_{i,\ell}$ in some $\rhos_i$
  whose endpoints $\rho_{i,j}^\pm$ and $\rho_{i,\ell}^\pm$ are nested,
  say $\rho_{i,j}^-<\rho_{i,\ell}^-<\rho_{i,\ell}^+<\rho_{i,j}^+$, to
  break apart and recombine into a pair of chords
  $\rho_a=(\rho_{i,j}^-,\rho_{i,\ell})^+$ and
  $\rho_b=(\rho_{i,\ell}^-,\rho_{i,j}^+)$. This gives an \emph{odd
    shuffle curve} at $e\infty$ and an element of
  $\cM^B(\x,\y;\vec{\rhos}')$ where $\vec{\rhos}'$ is
  obtained from $\vec{\rhos}$ by replacing $\rho_{i,j}$ and
  $\rho_{i,\ell}$ in $\rhos_i$ with $\rho_a$ and $\rho_b$.
\end{enumerate}
See Figure~\ref{fig:degens-first-three} for examples of the first
three kinds of degenerations.
The fourth is more difficult to
visualize; see~\cite{LOT1}.

\begin{remark}
  This analytic setup builds on the ``cylindrical reformulation''
  of Heegaard Floer theory~\cite{Lipshitz06:CylindricalHF}. It
  relates to the original formulation of Heegaard Floer theory, in
  terms of holomorphic disks in $\Sym^g(\Sigma)$, by thinking of a map
  $\bD\to\Sym^g(\Sigma)$ as a multi-valued map $\bD\to\Sigma$ and
  then taking the graph. See, for instance,
  \cite[Section~13]{Lipshitz06:CylindricalHF}. Some of the results were
  previously proved in~\cite{Lipshitz06:BorderedHF}.
\end{remark}
\subsection{Type D modules}\label{sec:type-D}
Fix a bordered Heegaard diagram $\HD$ and a suitable almost complex
structure~$J$.
Let $\PMC=-\bdy \HD$ be the
orientation reverse of the pointed matched circle given by
$\bdy\HD$. Given a generator $\x\in\Gen(\HD)$, let $I_D(\x)$ denote
the indecomposable idempotent of $\Alg(\PMC,0)\subset\Alg(\PMC)$
corresponding to the set of $\alpha$-arcs which are \emph{disjoint}
from~$\x$. This gives an action of the sub-ring $\Idem(\PMC)$ of
$\Alg(\PMC)$ on $X(\HD)$. As a module,
$\CFDa(\HD)=\Alg(\PMC)\otimes_{\Idem(\PMC)}X(\HD).$
Define a differential on $\CFDa(\HD)$ by 
\begin{multline*}
  \!\!\!\!\bdy \x=\sum_{\y\in\Gen(\HD)}\sum_{\begin{subarray}{c}B\in\pi_2(\x,\y)\\(\rho_1,\dots,\rho_n)\end{subarray}}\#\cM^B(\x,\y;(\{\rho_1\},\dots,\{\rho_n\}))\\[-16pt]
  \cdot a(\{-\rho_1\})\cdots
  a(\{-\rho_n\})\cdot \y,\!\!\!\!
\end{multline*}
with the convention that the number of elements in an
infinite set is zero. Here $-\rho$ denotes a chord $\rho$ in $(\bdy\HD,\CircPts)$ with
orientation reversed, so as to be a chord in $\PMC$. (To ensure
finiteness of these sums, we need to
impose an additional condition on the Heegaard diagram~$\HD$, called
\emph{provincial admissibility}~\cite{LOT1}.) Extend $\bdy$ to
all of $\CFDa(\HD)$ by the Leibniz rule.

\begin{theorem}[\cite{LOT1}]\label{thm:CFD-satisfies}
  The map $\bdy$ is a differential, i.e., $\bdy^2=0$.
\end{theorem}
\begin{proofsketch}
Let $a_{\x,\y}$ denote the coefficient of $\y$ in $\bdy(\x)$. The
equation $\bdy^2(\x)=0$ simplifies to the equation that, for all $\y$,
\[
d(a_{\x,\y})+\sum_\w a_{\x,\w}a_{\w,\y}=0.
\]
As usual in Floer theory, we prove this by considering
the boundary of the $1$-dimensional moduli spaces of curves. Of the four types of
degenerations, type~\ref{item:shuffle} does not
occur, since each $\rhos_i$ is a singleton set.
Type~\ref{item:two-story} gives the terms of the form
$a_{\x,\w}a_{\w,\y}$.  Type~\ref{item:split} with
a split curve degenerating corresponds to
$d(a_{\x,\y})$. Type~\ref{item:split} with no
split curves
cancel in pairs against themselves and
type~\ref{item:join}.  See also~\cite{LOT0}.
\end{proofsketch}

\begin{theorem}[\cite{LOT1}]\label{thm:D-invariance}
  Up to homotopy equivalence, the mod\-ule $\CFDa(\HD)$ is independent
  of the (provincially admissible) bordered Heegaard diagram $\HD$
  representing the bordered $3$-manifold $Y$.
\end{theorem}
The proof is similar to the closed case \cite[Theorem
6.1]{OS04:HolomorphicDisks}.
Theorem~\ref{thm:D-invariance} justifies writing $\CFDa(Y)$ for the
homotopy equivalence class of $\CFDa(\HD)$ for any bordered
Heegaard diagram $\HD$ for~$Y$.

\begin{remark}
  Since $\pi_2(\x,\y)$ is empty unless $\spinc(\x)=\spinc(\y)$,
  $\CFDa(Y)$ decomposes as a direct sum over $\SpinC$-structures on $Y$.
\end{remark}

\begin{remark}
  The modules $\CFDa(\HD)$ are projective over $\Alg(\PMC)$. Thus, the
  notions of homotopy equivalence and quasi-isomorphism are the same
  for these modules.
\end{remark}
\subsection{Type A modules}\label{sec:type-A}
Again, fix a bordered Heegaard diagram $\HD$ and a suitable almost
complex structure~$J$, but now let $\PMC=\bdy
\HD$.
Given a
generator $\x\in\Gen(\HD)$, let $I_A(\x)$ denote the indecomposable
idempotent in $\Alg(\PMC,0)\subset\Alg(\PMC)$ corresponding to the set
of $\alpha$-arcs intersecting $\x$ (opposite of
$I_D(\x)$), again making $X(\HD)$ into a module over
$\Idem(\PMC)\subset\Alg(\PMC)$.

Define an $\Ainf$ action of $\Alg(\PMC)$ on $X(\HD)$ by setting
\begin{multline*}
m_{n+1}(\x,a(\rhos_1),\dots,a(\rhos_n))\\=
\sum_{\y\in\Gen(\HD)}\sum_{B\in\pi_2(\x,\y)}\#\cM^B(\x,\y;(\rhos_1,\dots,\rhos_n))\cdot \y,
\end{multline*}
and extending multi-linearly. As for $\CFDa$, to ensure finiteness of
these sums, we need to assume that $\HD$ is provincially admissible.

\begin{theorem}\label{thm:CFA-satisfies}
  The operations $m_{n+1}$ satisfy the $\Ainf$ module relation.
\end{theorem}
\begin{proofsketch}
  Since $\Alg(\PMC)$ is a differential algebra, the $\Ainf$ relation
  for $\CFAa(Y)$ takes the form:
  \begin{equation}\label{eq:A-Ainf-rel}
    \begin{split}
      0 &=
      \sum_{i+j=n+2}m_i(m_j(\x,a_1,\dots,a_{j-1}),\dots,a_{n})\\
      &\quad+\sum_{\ell=1}^{n} m_{n+1}(\x,a_1,\dots,\partial a_\ell,\dots,a_{n})\\
      &\quad+\sum_{\ell=1}^{n-1} m_{n}(\x,a_1,\dots,a_\ell a_{\ell+1},\dots,a_{n}).
    \end{split}
  \end{equation}
  The first term in Equation~[\ref{eq:A-Ainf-rel}] corresponds to
  degenerations of type~\ref{item:two-story}. The second term
  corresponds to degenerations of
  types~\ref{item:join} and \ref{item:shuffle}, depending on whether
  one of the strands in the crossing being resolved is horizontal
  (\ref{item:join}) or not
  (\ref{item:shuffle}). The third term
  corresponds to degenerations of type~\ref{item:split}.
\end{proofsketch}

\begin{theorem}[\cite{LOT1}]
  Up to $\Ainf$ homotopy equivalence, the $\Ainf$ module $\CFAa(\HD)$
  is independent of the (provincially admissible) bordered Heegaard
  diagram $\HD$ representing the bordered $3$-manifold $Y$.
\end{theorem}
Again, the proof is similar to the closed case
\cite[Theorem 6.1]{OS04:HolomorphicDisks}.

\begin{remark}
  Like $\CFDa(Y)$, the module $\CFAa(Y)$ breaks up as a sum over
  $\SpinC$-structures on $Y$. 
  Further, as a chain complex, $\CFAa(Y)$
  breaks up as a sum over relative $\SpinC$ structures on $(Y,\bdy Y)$
  (but the module structure does not respect this splitting).
\end{remark}
\begin{remark}
  Over the algebra $\Alg(\PMC)$, $\Ainf$ homotopy equi\-val\-ence and
  $\Ainf$ quasi-isomorphism are equivalent notions. (This is true
  quite generally for $\Ainf$ modules over $\Ainf$ algebras over fields.)
\end{remark}
\begin{remark}
  It is always possible to choose a Heegaard diagram $\HD$ for $Y$ so
  that the higher products $m_n$, $n>2$, vanish on $\CFAa(\HD)$, so
  that $\CFAa(\HD)$ is an honest differential module. One way to do so
  is using an analogue of \emph{nice
    diagrams}~\cite{SarkarWang07:ComputingHFhat}.
\end{remark}

\subsection{Bimodules}
Next, suppose $Y$ is a \emph{strongly bordered $3$-manifold with two
boundary components}. By this we mean that we have a
$3$-manifold $Y$ with
boundary decomposed as $\bdy Y=\bdy_LY\amalg\bdy_RY$, homeomorphisms
$\phi_L\co F(\PMC_L)\to \bdy_LY$ and $\phi_R\co F(\PMC_R)\to \bdy_RY$,
and a framed arc $\gamma$ connecting the basepoints in $F(\PMC_L)$ and
$F(\PMC_R)$ and pointing into the preferred disks of $F(\PMC_L)$ and
$F(\PMC_R)$. Associated to $Y$ are bimodules
$\lsub{\Alg(-\PMC_L),\Alg(-\PMC_R)}\CFDDa(Y)$;
$\lsub{\Alg(-\PMC_L)}\CFDAa(Y)_{\Alg(\PMC_R)}$; and
$\CFAAa(Y)_{\Alg(\PMC_L),\Alg(\PMC_R)}$ defined by treating, respectively, both
$\bdy_LY$ and $\bdy_RY$ in type $D$ fashion; $\bdy_LY$ in type $D$
fashion and $\bdy_RY$ in type $A$ fashion; and both $\bdy_LY$ and
$\bdy_RY$ in type $A$ fashion.

An important special case of $3$-manifolds with two boundary
components is \emph{mapping cylinders}. Given an isotopy class of maps
$\phi\co F(\PMC_1)\to F(\PMC_2)$ taking the distinguished disk of
$F(\PMC_1)$ to the distinguished disk of $F(\PMC_2)$ and the basepoint
of $F(\PMC_1)$ to the basepoint of $F(\PMC_2)$---called a
\emph{strongly based mapping class}---the mapping cylinder $M_\phi$ of
$\phi$ is a strongly bordered $3$-manifold with two boundary
components. Let $\CFDDa(\phi)=\CFDDa(M_\phi)$,
$\CFDAa(\phi)=\CFDAa(M_\phi)$ and $\CFAAa(\phi)=\CFAAa(M_\phi)$.

The set of strongly based mapping classes forms a groupoid, with
objects the pointed matched circles representing genus $g$ surfaces
and $\Hom(\PMC_1,\PMC_2)$ the strongly based mapping
classes $F(\PMC_1)\to F(\PMC_2)$. In particular, the automorphisms of
a particular pointed matched circle $\PMC$ form the (strongly based)
mapping class group.

\begin{remark}
  The functors $\CFDAa(\phi)\DT\cdot$ give an action of the strongly
  based mapping class group on the (derived) category of left
  $\Alg(\PMC)$-modules; this action categorifies the standard action
  on $\Lambda^*H_1(F(\PMC);\Field)$.
  This action is faithful~\cite{LOT:faith}.
\end{remark}
\section{Gradings}

Let $\PMC=(Z,\CircPts,M,z)$ denote the genus $1$ pointed matched
circle from Example~\ref{eg:algebras}. Consier the following elements
of $\Alg(\PMC,1)$:
\begin{align*}
  x &= a(\{1,2\},\{2,3\},(1\mapsto 3,2\mapsto2))\\
  y &= a(\{1,2\},\{1,4\},(1\mapsto 1,2\mapsto4)).
\end{align*}
A short computation shows that $y\cdot x=d((dx)\cdot y)$.
It follows that there is no $\ZZ$-grading on $\Alg(\PMC,1)$ with homogeneous
basic generators.
A similar argument applies to $\Alg(\PMC',i)$ for any~$\PMC'$, as long
as $\Alg(\PMC',i)$
involves at least two moving strands.

There is, however, a grading in a more complicated sense. Let $G$ be a group and $\lambda\in G$ a distinguished central element. A grading of a differential algebra $\Alg$ by $(G,\lambda)$ is a decomposition 
$
\Alg=\bigoplus_{g\in G}\Alg_g
$
so that $d(\Alg_g)\subset \Alg_{\lambda^{-1} g}$ and $\Alg_g\cdot
\Alg_h\subset \Alg_{gh}$. Taking $G=\ZZ$ and $\lambda=1$ recovers the
usual notion of a $\ZZ$-grading of homological type. 

The corresponding notion for modules is a grading by a $G$-set. A
grading of a left differential $\Alg$-module $M$ by a left $G$-set $S$
is a decomposition
$
M=\bigoplus_{s\in S}M_s
$
so that $\Alg_gM_s\subset M_{gs}$ and $\bdy(M_s)\subset M_{\lambda^{-1} s}$.
Similarly, right modules are graded by right $G$-sets. If $M$ is
graded by a $G$-set $S$ and $N$ is graded by $T$ then $M\otimes_\Alg
N$ is graded by $S\times_G T$, which retains an action of~$\lambda$.
These more
general kinds of gradings have been considered by, e.g.,
N\u ast\u asescu and Van Oystaeyen~\cite{GroupGradings}.

Given a surface $F$, let $G(F)$ be the $\ZZ$-central extension of
$H_1(F)$, 
\[
\ZZ \overset{\lambda}{\hookrightarrow} \smallGroup(F)
  {\twoheadrightarrow} H_1(F),
\]
where $1 \in \ZZ$ maps to $\lambda \in
\smallGroup(\PtdMatchCirc)$. Explicitly, $G(F)$ is a subgroup of the
group $\OneHalf\ZZ\times H_1(F)$ with multiplication
\[
  (k_1, \alpha_1) \cdot (k_2, \alpha_2) = (k_1 + k_2 +
  \alpha_1 \cap \alpha_2, \alpha_1 + \alpha_2).
\]
Here, $\cap$ denotes the intersection pairing on $H_1(F)$.
It turns out that $\Alg(\PMC)$ has a grading
by $(G(F(\PMC)),\lambda)$~\cite[Section~3.3]{LOT1}.
Similarly, given a $3$-manifold $Y$ bordered by $F(\PMC)$, one can
construct $G(F(\PMC))$-set gradings on $\CFAa(Y)$ and
$\CFDa(-Y)$~\cite{LOT1}.

  Even in the closed case, the grading on Heegaard Floer homology has
  a somewhat nonstandard form: a \emph{partial relative cyclic grading}. That
  is, generators do not have well-defined gradings, 
  but only
  well-defined grading differences $\gr(\x,\y)$; the grading
  difference $\gr(\x,\y)$ is only
  defined for 
  generators representing the
  same $\SpinC$-structure; and $\gr(\x,\y)$ is well-defined only
  modulo 
  the divisibility of $c_1(\spinc(\x))$.
  A partial
  relative cyclic grading is precisely a grading by a $\ZZ$-set. This
  leads naturally to a graded version of the pairing
  theorems, including Equation~[\ref{eq:pairing-1}]~\cite{LOT1}.

\section{Deforming the diagonal and the pairing theorems}\label{sec:pairing}
The tensor product pairing theorems are the main motivation for the
definitions of the modules and bimodules. 
We will sketch the proof of the
archetype, Equation~[\ref{eq:pairing-1}].
Fix bordered Heegaard diagrams $\HD_1$ and $\HD_2$ for $Y_1$ and $Y_2$,
respectively, with $\bdy\HD_1=-\bdy\HD_2$. It is easy to see that $\HD=\HD_1\cup_\bdy\HD_2$ is a
Heegaard diagram for $Y_1\cup_\bdy Y_2$.

There are two sides to the proof, one algebraic and one analytic. We
start with the algebra. Typically, the $\Ainf$ tensor product $M\DTP
N$ of $\Ainf$ modules $M$ and $N$ is defined using a chain complex
whose underlying vector space is $M\otimes_\Ground
T^*\!\Alg\otimes_\Ground N$ (where $T^*\!\Alg$ is the tensor algebra of
$\Alg$, and $\Ground$ is the ground ring of $\Alg$---for us, the ring
of idempotents). This complex is typically infinite-dimensional, and
so unlikely to align easily with $\CFa$.

However, the differential module $\CFDa(\HD_2)$ has a special form: it
is given as $\Alg(\PMC)\otimes_{\Idem(\PMC)}X(\HD_2)$, so the
differential is encoded by a map
$
\delta^1\co X(\HD_2)\to \Alg(\PMC)\otimes_{\Idem(\PMC)}X(\HD_2).
$
This allows us to define a smaller model for the $\Ainf$ tensor
product. Let $\delta^n\co X(\HD_2)\to \Alg(\PMC)^{\otimes
  n}\otimes X(\HD_2)$ be the result of iterating $\delta^1$
$n$ times.
For notational convenience, let $M=\CFAa(\HD_1)$ and
$X=X(\HD_2)$.
Define $\CFAa(\HD_1)\DT\CFDa(\HD_2)$ to be the $\Field$-vector
space $M\otimes_{\Idem(\PMC)}X$, with differential (graphically
depicted in
Figure~\ref{fig:DT})
\begin{equation}
  \bdy=\sum_{i=0}^\infty (m_{i+1}\otimes \Id_X)\circ
  (\Id_M\otimes\delta^i).
\label{eq:def-DT}
\end{equation}

\begin{figure}
  \centering
  $\begin{tikzpicture}
    \node at (0,0) (t) {};
    \node at (0,-1.5) (delta) {$\delta^1$};
    \node at (0,-3) (b) {};
    \node at (-.5,-3) (bl) {};
    \draw[Dmodar] (t) to (delta);
    \draw[Dmodar] (delta) to (b);
    \draw[algarrow] (delta) to (bl);
    \node at (-.25, -3.5) (letter) {(a)};
  \end{tikzpicture}$\qquad\qquad
  $\begin{tikzpicture}
    \node at (0,0) (t) {};
    \node at (.5,0) (t2) {};
    \node at (1,0) (t3) {};
    \node at (1.5, 0) (t4) {};
    \node at (0,-1.5) (m) {$m_4$};
    \node at (0,-3) (b) {};
    \draw[Amodar] (t) to (m);
    \draw[Amodar] (m) to (b);
    \draw[algarrow] (t2) to (m);
    \draw[algarrow] (t3) to (m);
    \draw[algarrow] (t4) to (m);
    \node at (.75, -3.5) (letter) {(b)};
  \end{tikzpicture}$\qquad
  $\begin{tikzpicture}
    \node at (0,0) (tlA) {};
    \node at (.5,0) (trA) {};
    \node at (0, -2) (mA) {$m_1$};
    \node at (0,-3) (blA) {};
    \node at (.5,-3) (brA) {};
    \draw[Dmodar] (trA) to (brA);
    \draw[Amodar] (tlA) to (mA);
    \draw[Amodar] (mA) to (blA);
    \node at (1,-1.5) (plus1) {+};
    \node at (1.5,0) (tlB) {};
    \node at (2,0) (trB) {};
    \node at (2, -1) (delta1B) {$\delta^1$};
    \node at (1.5, -2) (mB) {$m_2$};
    \node at (1.5,-3) (blB) {};
    \node at (2,-3) (brB) {};
    \draw[Dmodar] (trB) to (delta1B);
    \draw[Dmodar] (delta1B) to (brB);
    \draw[Amodar] (tlB) to (mB);
    \draw[Amodar] (mB) to (blB);
    \draw[algarrow] (delta1B) to (mB);
    \node at (1.75, -3.5) (letter) {(c)};
    \node at (2.5,-1.5) (plus2) {+};
    \node at (3,0) (tlC) {};
    \node at (3.5,0) (trC) {};
    \node at (3.5, -.66) (delta1C) {$\delta^1$};
    \node at (3.5, -1.33) (delta2C) {$\delta^1$};
    \node at (3, -2) (mC) {$m_3$};
    \node at (3,-3) (blC) {};
    \node at (3.5,-3) (brC) {};
    \draw[Dmodar] (trC) to (delta1C);
    \draw[Dmodar] (delta1C) to (delta2C);
    \draw[Dmodar] (delta2C) to (brC);
    \draw[Amodar] (tlC) to (mC);
    \draw[Amodar] (mC) to (blC);
    \draw[algarrow] (delta1C) to (mC);
    \draw[algarrow] (delta2C) to (mC);
    \node at (4,-1.5) (plus3) {$+\cdots$};
  \end{tikzpicture}
  $
  \caption{\textbf{The operation $\DT$.} (a) Graphical representation
    of $\delta^1$ for $\CFDa$. (b) Graphical representation of
    $m_{n+1}$ ($n=3$) for $\CFAa$. (c) Graphical representation of the
    sum from Equation~[\ref{eq:def-DT}]. In all cases, dashed lines
    represent module elements and solid lines represent algebra
    elements.}
  \label{fig:DT}
\end{figure}

The sum in Equation~[\ref{eq:def-DT}] is not \emph{a priori}
finite. To ensure that it is finite, we need to assume an additional
boundedness condition on either $\CFAa(\HD_1)$ or
$\CFDa(\HD_2)$. These boundedness conditions correspond to an
admissibility hypothesis for $\HD_i$, which in turn guarantees that
$\HD_1\cup_\bdy\HD_2$ is (weakly) admissible.

\begin{lemma}
  There is a canonical homotopy equivalence
  \[
  \CFAa(\HD_1)\DTP_{\Alg(\PMC)}\CFDa(\HD_2)\simeq
  \CFAa(\HD_1)\DT_{\Alg(\PMC)}\CFDa(\HD_2).
  \]
\end{lemma}
The proof is straightforward.

We turn to the analytic side of the argument next.  Because of how the
idempotents act on $\CFAa(\HD_1)$ and $\CFDa(\HD_2)$, there is an
obvious identification between generators $\x_L\otimes\x_R$ of
$\CFAa(\HD_1)\DT\CFDa(\HD_2)$ and generators $\x$ of $\CFa(\HD)$.  

Let $Z=\bdy\HD_1=\bdy\HD_2\subset\HD$. The differential on $\CFa(\HD)$
counts rigid $J$-holomorphic curves in
$\Sigma\times[0,1]\times\RR$. For a sequence of almost complex
structures $J_r$ with longer and longer necks at $Z$, such curves
degenerate to pairs of curves $(u_L,u_R)$ for $\HD_1$ and $\HD_2$,
with matching asymptotics at $e\infty$. More precisely, in the
limit
as we stretch the neck, the moduli space degenerates to a fibered
product
\begin{multline}\label{eq:fiber-prod}
\coprod_{\rho_1,\dots,\rho_i}\!\!\!\!\cM(\x_L,\y_L,(\rho_1,\dots,\rho_i))\\[-8pt]
\sos{\ev_L}{\times}{\ev_R}\cM(\x_R,\y_R,(\rho_1,\dots,\rho_i)).
\end{multline}
Here, 
$
\ev_L\co \cM(\x,\y;(\rho_1,\dots,\rho_i))\to\RR^i/\RR
$
is given by taking the successive height differences (in the
$\RR$-coordinate) of the chords $\rho_1,\dots,\rho_i$, and similarly
for $\ev_R$. Also, we are suppressing homology classes from the
notation.

Since we are taking a fiber product over a large\hyp dimensional space,
the moduli spaces in Equation~[\ref{eq:fiber-prod}] are not
conducive to defining invariants of $\HD_1$ and $\HD_2$. To deal with
this, we deform the matching condition, considering instead the fiber
products
\begin{multline*}
\coprod_{\rho_1,\dots,\rho_i}\!\!\!\!\cM(\x_L,\y_L,(\rho_1,\dots,\rho_i))\\[-8pt]
  \sos{R\cdot\ev_L}{\times}{\ev_R}\cM(\x_R,\y_R,(\rho_1,\dots,\rho_i))
\end{multline*}
and sending $R\to\infty$. In the limit, some of the chords on the left
collide, while some of the chords on the right become infinitely far
apart. The result exactly recaptures the definitions of $\CFAa(\HD_L)$
and $\CFDa(\HD_R)$ and the algebra of Equation~\eqref{eq:def-DT}~\cite{LOT1}.


\begin{remark}\label{rmk:nice-pairing}
  In~\cite{LOT1}, we also give another proof of the pairing
  theorem~[\ref{eq:pairing-1}], using the technique of nice
  diagrams~\cite{SarkarWang07:ComputingHFhat}.
\end{remark}
\section{Computing with bordered Floer homology}\label{sec:computing}
\subsection{Computing $\HFa$}
\begin{figure}
  \centering
  \input{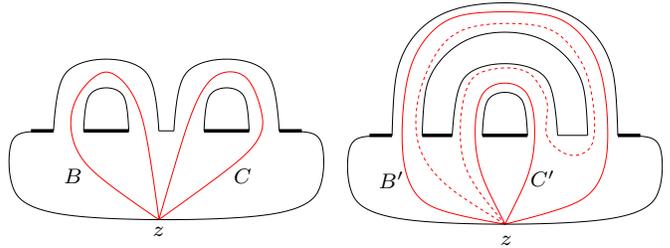}
  \caption{\textbf{A local picture of an arc-slide diffeomorphism.} Left: a
pair of pants with two distinguished curves $B$ and $C$. Right:
another pair of pants with
distinguished curves $B'$ and $C'$. The arc-slide diffeomorphism
carries $B$ to the dotted curve on the right and
$C$ on the left to $C'$ on the right.
This diffeomorphism can be extended to a diffeomorphism between
surfaces associated to pointed matched circles: in such a surface
there may be further handles attached along the four dark
intervals.
This figure also appears in~\cite{LOT4}.}
  \label{fig:arcslides}
\end{figure}

Let $Y$ be a closed $3$-manifold. As discussed earlier, $Y$ admits a
Heegaard splitting into two handlebodies, glued by some homeomorphism
$\phi$ between their boundaries. Via the pairing
theorems~[\ref{eq:MorPair1}] and~[\ref{eq:MorPair2}], this reduces
computing $\HFa(Y)$ to computing
$\CFDa(\HB_g)$ for some particular bordered
  handlebody $\HB_g$ of each genus $g$ and
$\CFDDa(\phi)$ for arbitrary $\phi$ in the
  strongly based mapping class group. 
For a appropriate~$\HB_g$, $\CFDa(\HB_g)$ is easy to compute.
Moreover, by~[\ref{eq:MorPair2}]
we
do not need to compute $\CFDDa(\phi)$ for every mapping class, just
for generators for the mapping class groupoid.  This
groupoid has a natural set of generators: arc-slides 
(see
Figure~\ref{fig:arcslides}, and
compare~\cite{Bene08:ChordDiagrams,Bene09:MCFactorNielsen}). 

It turns out that the type \DD\
invariants of arc-slides are determined by a small amount of geometric
input (essentially, the set of generators and a non-degeneracy
condition for the differential) and the condition that $\bdy^2=0$~\cite{LOT4}.

These techniques also allow one to compute all types of the bordered
invariants for any bordered $3$-manifold.

\subsection{Cobordism maps}\label{sec:cobordism-maps}
Next, we discuss how to compute the map $\hat{f}_W\co \CFa(Y_1)\to
\CFa(Y_4)$ associated to a 4-dimensional cobordism $W$ from $Y_1$ to
$Y_4$. The cobordism $W$ can be decomposed into three
cobordisms $W_1W_2W_3$ where $W_i\co Y_i\to Y_{i+1}$ consists 
of $i$-handle attachments and $\hat{f}_{W}$ is a
corresponding composition $\hat{f}_{W_3}\circ
\hat{f}_{W_2}\circ \hat{f}_{W_1}$. 

The maps $\hat{f}_{W_1}$ and $\hat{f}_{W_3}$ are simple to describe:
$Y_2\cong Y_1\#^k(S^2\times S^1)$, while $Y_3\cong Y_4\#^\ell
(S^2\times S^1)$;  $\CFa(S^2\times S^1)$ is (homotopy equivalent to)
$\Field\oplus\Field=H_*(S^1;\Field)$; and the invariant $\HFa(Y)$
satisfies a K\"unneth theorem for connect sums, so
\[
\CFa(Y_1\#^k(S^2\times S^1))\cong \CFa(Y_1)\otimes_\Field H_*(T^k;\Field)
\]
(with respect to appropriate Heegaard diagrams),
where $T^k = (S^1)^k$ is the $k$-dimensional torus.
Let $\theta$ be the top-dimensional generator of $H_*(T^k;\Field)$ and
$\eta$ the bottom-dimensional generator of $H_*(T^\ell;\Field)$. Then
$\hat{f}_{W_1}$ is $\x\mapsto \x \otimes\theta$ while
$\hat{f}_{W_3}$ takes $\x\otimes \eta\mapsto \x$ and
$\x\otimes\epsilon\mapsto 0$ if $\gr(\epsilon)>\gr(\eta)$.

By contrast, $\hat{f}_{W_2}$ is defined by counting
holomorphic triangles in a suitable Heegaard triple-diagram. 
%
Two additional properties of
bordered Floer theory allow us to compute $\hat{f}_{W_i}$:
\begin{itemize}
\item The invariant $\CFDa(\HB)$ of a handlebody $\HB$ is rigid, in
  the sense that it has no nontrivial graded automorphisms. 
  This allows one to compute the homotopy equivalences between the
  results of making different choices in the computation of $\CFa(Y)$.
\item There is a pairing theorem for holomorphic triangles.
\end{itemize}
Given these, one can compute $\hat{f}_W$ as follows. Using results
from the previous section, we can compute $\CFa(Y_2)$ (respectively
$\CFa(Y_3)$) using a Heegaard decomposition making the decomposition
$Y_2\cong Y_1\#^k(S^2\times S^1)$ (respectively $Y_3\cong
Y_4\#^\ell(S^2\times S^1)$) manifest. With respect to this
decompositions, the map $\hat{f}_{W_1}$ (respectively $\hat{f}_{W_3}$)
is easy to read off.

To compute $\hat{f}_{W_2}$ one works with 
Heegaard decompositions of $Y_2$ and $Y_3$ with respect to which the
cobordism $W_2$ takes a particularly simple form, replacing one of the
handlebodies $\HB$ of a Heegaard decomposition of~$Y_2$ with a
differently framed handlebody $\HB'$. It is
easy to compute the triangle map $\CFDa(\HB)\to \CFDa(\HB')$. By the
pairing theorem for triangles, extending this
map by the identity
map on the rest of the decomposition gives the map $\hat{f}_{W_2}$.
Finally, the rigidity result allows one to write down the isomorphisms
between $\CFa(Y_2)$ (and $\CFa(Y_3)$) computed in the two
different ways. The map $\hat{f}_W$ is then the composition of the
three maps $\hat{f}_{W_i}$ and the equivalences connecting
the two different models of $\CFa(Y_2)$ and of $\CFa(Y_3)$.

For more details, the reader is referred to~\cite{LOTCobordisms}.

\subsection{Polygon maps and the Ozsv\'ath-Szab\'o spectral sequence}
Mikhail Khovanov introduced a
categorification of the Jones
polynomial~\cite{Khovanov00:CatJones}. This categorification
associates to an oriented link $L$ a bigraded abelian group $\Kh_{i,j}(L)$, the
\emph{Khovanov homology of $L$}, whose graded Euler characteristic  is $(q+q^{-1})$ times the Jones polynomial $J(L)$.
There is also a reduced
version $\rKh(L)$, whose graded Euler characteristic is simply
$J(L)$. 
(The reduced theory requires that one mark one component of
$L$.) 
The skein relation for $J(L)$ is replaced by a \emph{skein exact
  sequence}. Given a link $L$ and a crossing $c$ of $L$, let $L_0$ and
$L_1$ be the two resolutions of~$c$:
\[
\includegraphics[scale=.666666]{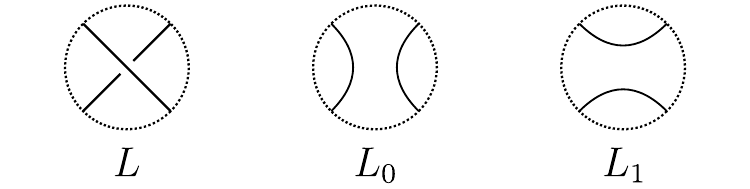}. 
\]
Then there is a long exact sequence relating the (reduced) Khovanov
homology groups of $L$, $L_1$, and $L_0$.

Szab\'o and the second author observed that the
Heegaard Floer group $\HFa(D(L))$ of the double cover of $S^3$
branched over $L$ satisfies a similar skein exact triangle to
(reduced) Khovanov homology, and takes the same value on an
$n$-component unlink (with some collapse of gradings). Using
these observations, they produced a spectral sequence from Khovanov
homology (with $\Field$-coefficients) to $\HFa(D(L))$~\cite{BrDCov}.
(Because of a difference in conventions, one must take the Khovanov
homology of the mirror $r(L)$ of $L$.)
J.~Baldwin recently showed~\cite{Baldwin:ss} that the entire spectral sequence
$\rKh(r(L)) \rightrightarrows \HFa(D(L))$
is a knot invariant.

Bordered Floer homology
can be used to compute this spectral sequence~\cite{LOT:DCov1,
  LOT:DCov2}. Write $L$ as the plat closure of
some braid $B$, and decompose $B$ as a product of braid generators
$s_{i_1}^{\epsilon_1}\cdots s_{i_k}^{\epsilon_k}$. The branched double
cover of a braid generator $s_i^\pm$ is the mapping cylinder of a Dehn
twist, and the branched double covers of the plats closing $B$ is a
handlebody~$\HB$. So 
$\CFa(D(K))$ is
quasi-isomorphic to
\begin{equation} \CFAa(\HB)\DT\CFDAa(\tau_{i_1})\DT\cdots\DT\CFDAa(\tau_{i_k})\DT\CFDa(\HB).\label{eq:cube-1}
\end{equation}

The bordered invariant of a Dehn twist $\tau_\gamma$ along
$\gamma\subset F(\PMC)$ can be written as a mapping cone of a map
between the identity cobordism $\Id=[0,1]\times F(\PMC)$ and
the manifold $Y_{0(\gamma)}$ obtained by $0$-surgery on $[0,1]\times F(\PMC)$ along~$\gamma$:
\begin{align}
  \CFDAa(\tau_\gamma)&\simeq \Cone(\CFDAa(\Id_{0(\gamma)})\to\CFDAa(\Id))\label{eq:cube-2}\\
  \CFDAa(\tau_\gamma^{-1})&\simeq \Cone(\CFDAa(\Id)\to \CFDAa(\Id_{0(\gamma)})).\label{eq:cube-3}
\end{align}
Applying  this observation to the tensor product in Formula~[\ref{eq:cube-1}] endows $\CFa(D(K))$ with 
with a filtration by $\{0,1\}^n$. The resulting spectral sequence has
$E_1$-page the Khovanov chain complex and $E_\infty$-page
$\HFa(D(K))$.

The next step is to compute the groups $\CFDAa(\Id_{0(\gamma)})$ for
the curves $\gamma$ corresponding to the braid generators and the maps
from [\ref{eq:cube-2}] and~[\ref{eq:cube-3}], which again
requires a small amount of geometry~\cite{LOT:DCov1}.

Finally, we identify this spectral sequence with the earlier spectral
sequence~\cite{BrDCov}. The key ingredient is another
pairing theorem identifying the algebra of tensor products of
mapping cones with counts of holomorphic polygons~\cite{LOT:DCov2}.





\section{Acknowledgements}
  We thank MSRI for hosting us in the spring of
  2010, during which part of this research was conducted.

  This is a slightly expanded and reformatted version of a paper which
  first appeared in the Proceedings of the National Academy of
  Sciences~\cite{LOT:tour}, available online at
  \url{http://www.pnas.org/content/early/2011/04/25/1019060108.short}.

\DeclareUrlCommand\eprint{\urlstyle{rm}}
\bibliographystyle{mypnas}
\bibliography{heegaardfloer}

\end{document}